\documentclass[11pt]{article}
\usepackage{graphicx}
\usepackage{ifpdf}    
\usepackage[letterpaper, hmargin=1in, vmargin=1.25in]{geometry}

    \renewcommand{\baselinestretch}{1.3}
  \renewcommand{\arraystretch}{1.1}
\begin{document}

 \title{A Heuristic Description of  Fast Fourier Transform}

 \author{Zhengjun Cao\,$^*$, \quad  Xiao Fan\\
  Department of Mathematics, Shanghai University, Shanghai,
  China.\\
  \textsf{$^*$\ \textsf{caozhj@shu.edu.cn}}  }

 \date{}\maketitle

\begin{abstract}

Fast Fourier Transform (FFT) is an efficient algorithm to compute the Discrete Fourier Transform (DFT) and its inverse.
 In this paper, we pay special attention to the description of  complex-data FFT. We analyze two common descriptions of FFT and  propose a new
presentation. Our heuristic description is helpful for  students and programmers to grasp the algorithm entirely and deeply.

\end{abstract}

\section{Introduction}

 Denote the set of all complex numbers by $\mathcal {C}$. Let $\omega=e^{2 \pi i/n}$ be a primitive $n$th root of unity. Let
  $$f(x)=\sum_{0\leq i<n}f_i\,x^i \in \mathcal {C}[x]  \eqno(1) $$
   be a polynomial of degree less than $n$ with its coefficient vector $(f_0, \cdots, f_{n-1})\in \mathcal {C}^n$. The $\mathcal {C}$-linear map
  $$
\mbox{\textbf{DFT}}_{\omega}:\left\{\begin{array}{lll}
  \mathcal {C}^n &\rightarrow &  \mathcal {C}^n\\
(f_0, \cdots, f_{n-1}) &\mapsto & \left( f(1), f(\omega), f(\omega^2),\cdots, f(\omega^{n-1})   \right)
\end{array}
 \right.
$$
which evaluates a polynomial at the powers of $\omega$ is called the Discrete Fourier Transform. Apparently,
the Discrete Fourier Transform is a special multipoint evaluation at the powers $1, \omega,  \omega^2, \cdots, \omega^{n-1}$.

A beginner will often evaluate the polynomial (1) in a manner that: first $f_{n-1}x^{n-1} $ is calculated, then
$f_{n-2}x^{n-2}, \cdots, f_1x$, and finally all of the terms of (1) are added together.  Such a process involves lots of multiplications and additions.
There is an elegant way to evaluate a polynomial, which is called Horner's rule. It evaluates $f(x)$ as follows: Rearrange this computation as
$$f(x)=(\cdots (f_{n-1}x+f_{n-2}\,)x+\cdots)\,x+f_0 \eqno(2) $$
Then start with $f_{n-1}$, multiply by $x$, add $f_{n-2}$, multiply by $x$, $\cdots$, multiply by $x$, add $f_0$. If we apply the Horner's rule to compute
$\mbox{\textbf{DFT}}_{\omega}$, it needs to run the evaluation algorithm $n$ times. Can we find more efficient algorithm to compute Discrete Fourier transform?

Fast Fourier Transform (FFT) is an efficient algorithm to compute the Discrete Fourier Transform (DFT) and its inverse.
It was popularized by a publication of J. Cooley and J. Tukey \cite{CT65} in 1965. It has been called the most important numerical algorithm of our lifetime.
 In this paper,we will pay special attention to the description of complex-data FFT. We will analyze two common descriptions of  FFT \cite{GG03,K97} and
 propose an  explicit and heuristic presentation of FFT.

\section{Analysis of two common descriptions of FFT}

\subsection{Description-1 }

We refer to Ref.[3] for
the first description of FFT.

Let $n=2^k \in \mathcal {N}$ with $k\in \mathcal {N}$,   $\omega \in \mathcal {C}$ be a primitive $n$th root of unity, and $f(x)=\sum_{0\leq i<n}f_i\,x^i\in \mathcal {C}[x]$ of degree less than $n$. To evaluate $f(x)$ at the powers $1, \omega, \omega^2, \cdots, \omega^{n-1}$, we divide $f(x)$ by $x^{n/2}-1$ and $x^{n/2}+1$ with remainder:
$$ f(x)=q_0(x)(x^{n/2}-1)+r_0(x)=q_1(x)(x^{n/2}+1)+r_1(x) \eqno(3)$$
for some $q_0(x), r_0(x), q_1(x), r_1(x)\in \mathcal {C}[x]$ of degree less than $n/2$.

Plugging in a power of $\omega$ for $x$ in (3), we find
\begin{eqnarray*}
f(\omega^{2\ell})&=&q_0(\omega^{2\ell})(\omega^{n\ell}-1)+r_0(\omega^{2\ell})=r_0(\omega^{2\ell})\\
f(\omega^{2\ell+1})&=&q_1(\omega^{2\ell+1})(\omega^{n\ell}\omega^{n/2}+1)+r_1(\omega^{2\ell+1})=r_1(\omega^{2\ell+1})
\end{eqnarray*}
for all $0\leq \ell < n/2$. We have used the facts that $\omega^{n}=1$ and $\omega^{n/2}=-1$. It  remains to evaluate $r_0(x)$ at the even powers of $\omega$ and $r_1(x)$ at the odd powers. Now $\omega^2$ is a primitive ($n/2$)th root of unity, and hence the first task is a DFT of order $n/2$. But also the evaluation of $r_1(x)$ can be reduced to a DFT of order $n/2$ by noting that $r_1(\omega^{2\ell+1})=r_1^*(\omega^{2\ell})$ for
$r_1^*(x)=r_1(\omega x)$. Since $n$ is a power of 2, we can proceed recursively to evaluate
$r_0(x)$ and $r_1^*(x)$ at the powers $1, \omega^2, \cdots, \omega^{2n-2}$ of $\omega^2$, and obtain the following algorithm.

\begin{table}[!h]
\tabcolsep 0pt
\caption{Description-1}
\vspace*{-12pt}
\begin{center}
\def\temptablewidth{1\textwidth}
{\rule{\temptablewidth}{.6pt}}

\begin{tabular*}{\temptablewidth}{@{\extracolsep{\fill}}l}
Input: $n=2^k \in \mathcal {N}$ with $k\in \mathcal {N}, f(x)=\sum_{0\leq i<n}f_i\,x^i\in \mathcal {C}[x]$,
and \\  \hspace*{8mm} the powers $ \omega, \omega^2, \cdots, \omega^{n-1}$ of a primitive $n$th root of unity $\omega \in \mathcal {C}$.\\
Output:  $\mathrm{DFT}_\omega(f)=(f(1),f(\omega),\cdots,f(\omega^{n-1}))\in \mathcal {C}^n.$\\
\hspace*{8mm}  1.  \textbf{if} $n=1$ \textbf{then return} $f_0$\\
  \hspace*{8mm}   2. $r_0(x)\longleftarrow \sum_{0\leq j<{n/2}}(f_j+f_{j+n/2})x^j$,\quad $r_1^*(x)\longleftarrow \sum_{0\leq j<{n/2}}(f_j-f_{j+n/2})\omega^j\,x^j $\\
\hspace*{8mm}  3. \textbf{call} the algorithm recursively to evaluate $r_0(x)$ and $r_1^*(x)$ at the powers of $\omega^2$\\
\hspace*{8mm}  4. \textbf{return} $\left(r_0(1), r_1^*(1), r_0(\omega^2), r_1^*(\omega^2), \cdots, r_0(\omega^{n-2}), r_1^*(\omega^{n-2})\right)$
\end{tabular*}
       {\rule{\temptablewidth}{.6pt}}
       \end{center}
       \end{table}

\subsection{Analysis of Description-1 }

To investigate the working flow in description-1, we take the following polynomial as an example:
$$f(x)=f_0+f_1x+f_2x^2+\cdots+f_7x^7 \in \mathcal {C}[x], \ \ \omega=e^{{2\pi i}/{8}}  \eqno(4) $$

At first,  we  have  the following two remainders:
\begin{eqnarray*}
 r_0(x)&=& (f_0+f_4)+(f_1+f_5)x+(f_2+f_6)x^2+(f_3+f_7)x^3  \hspace*{48mm} (5) \\
r_1^*(x)&=& (f_0-f_4)+(f_1-f_5)\omega\,x+(f_2-f_6)\omega^2\,x^2+(f_3-f_7)\omega^3\,x^3 \hspace*{35mm} (6)
\end{eqnarray*}
If we evaluate them at $(\omega^2)^0$, then $r_0(1)=f(1) $, $r_1^*(1)=f(\omega)$.

We now proceed to the step 3.   By (5)  we have the following two remainders:
\begin{eqnarray*}
& & [(f_0+f_4)+(f_2+f_6)]+[(f_1+f_5)+(f_3+f_7)]x  \hspace*{58mm} (7)\\
& & [(f_0+f_4)-(f_2+f_6)]+[(f_1+f_5)-(f_3+f_7)]\omega\,x \hspace*{55mm} (8)
\end{eqnarray*}
Likewise, by (6) we have the following two remainders:
\begin{eqnarray*}
& &  [(f_0-f_4)+(f_2-f_6)\omega^2]+[(f_1-f_5)\omega+(f_3-f_7)\omega^3] x \hspace*{48mm} (9)\\
& &  [(f_0-f_4)-(f_2-f_6)\omega^2]+[(f_1-f_5)\omega-(f_3-f_7)\omega^3]\omega\, x \hspace*{44mm} (10)
\end{eqnarray*}
Should we evaluate them at $(\omega^2)^1$? If that we can not find correct answers.

\emph{Drawback}: The sentence that, ``call the algorithm recursively to evaluate $r_0(x)$ and $r_1^*(x)$ at the powers of $\omega^2$",
is really too vague to specify the working flow.

\subsection{Description-2 }
We refer to Ref.[4] for
the second description of FFT.

Let $n=2^k \in \mathcal {N}$ with $k\in \mathcal {N}$,   $\omega \in \mathcal {C}$ be a primitive $n$th root of unity, and $f(x)=\sum_{0\leq i<n}f_i\,x^i\in \mathcal {C}[x]$ of degree less than $n$. To evaluate $f(x)$ at the powers $1, \omega, \omega^2, \cdots, \omega^{n-1}$, we carry out the following scheme. (In these formulas the parameters $s_j$ and $t_j$ are either 0 or 1, so that each ``pass" represents $2^k$ elementary computations.)

       \begin{table}[!h]
\tabcolsep 0pt
\caption{Description-2}
\vspace*{-12pt}
\begin{center}
\def\temptablewidth{1\textwidth}
{\rule{\temptablewidth}{.6pt}}

\begin{tabular*}{\temptablewidth}{@{\extracolsep{\fill}}l}
\textbf{Pass 0.} Let $A^{[0]}(t_{k-1},\cdots,t_0)=f_t,$ where $t=(t_{k-1},\cdots,t_0)_2$.\\
\textbf{Pass 1.} Set $A^{[1]}(s_{k-1},t_{k-2},\cdots,t_0)\leftarrow $\\
\hspace*{18mm}$A^{[0]}(0,t_{k-2},\cdots,t_0)+\omega^{2^{k-1}s_{k-1}}A^{[0]}(1,t_{k-2},\cdots,t_0)$\\
\textbf{Pass 2.} Set $A^{[2]}(s_{k-1},s_{k-2},t_{k-3},\cdots,t_0)\leftarrow $\\
\hspace*{18mm} $ A^{[1]}(s_{k-1},0,t_{k-3},\cdots,t_0)+\omega^{2^{k-2}(s_{k-2}s_{k-1})_2}A^{[1]}(s_{k-1},1,t_{k-3},\cdots,t_0)$\\
\hspace*{18mm} $\cdots$\\
\textbf{Pass k.} Set $A^{[k]}(s_{k-1},\cdots,s_1,s_0)\leftarrow$\\
\hspace*{18mm} $ A^{[k-1]}(s_{k-1},\cdots,s_1,0)+\omega^{(s_0s_1\cdots s_{k-1})_2}A^{[k-1]}(s_{k-1},\cdots,s_1,1)$
\end{tabular*}
       {\rule{\temptablewidth}{.6pt}}
       \end{center}
       \end{table} \vspace*{-8mm}

       It is  easy to find that
$$A^{[k]}(s_{k-1},\cdots,s_1,s_0)=f(\omega^s),\quad  \mathrm{where}\ \  s=(s_0,s_1,\cdots,s_{k-1})_2 \eqno(11)$$
Notice that the binary digits of $s$ are reversed in the final result (11).

   \subsection{Analysis of Description-2 }

 To investigate the working flow in description-2,  we also take the polynomial (4) as an example.
 Clearly, we have

 { \renewcommand{\baselinestretch}{.5}
  \renewcommand{\arraystretch}{.3}  \normalsize \small
 \begin{eqnarray*}
& &  A^{[0]}(000)=f_0,\ \ A^{[0]}(001)=f_1,\ \  A^{[0]}(010)=f_2,\ \  A^{[0]}(011)=f_3\ \\
& &  A^{[0]}(100)=f_4,\ \  A^{[0]}(101)=f_5,\ \  A^{[0]}(110)=f_6,\ \ A^{[0]}(111)=f_7\
 \end{eqnarray*}
 In pass 1,
 \begin{eqnarray*}
& &  A^{[1]}(000)=A^{[0]}(000)+A^{[0]}(100)=f_0+f_4,\ \  A^{[1]}(100)=A^{[0]}(000)+\omega^{2^2}A^{[0]}(100)=f_0-f_4 \\
& &  A^{[1]}(001)=A^{[0]}(001)+A^{[0]}(101)=f_1+f_5,\ \  A^{[1]}(101)=A^{[0]}(001)+\omega^{2^2}A^{[0]}(101)=f_1-f_5 \\
& &  A^{[1]}(010)=A^{[0]}(010)+A^{[0]}(110)=f_2+f_6,\ \  A^{[1]}(110)=A^{[0]}(010)+\omega^{2^2}A^{[0]}(110)=f_2-f_6 \\
& &  A^{[1]}(011)=A^{[0]}(011)+A^{[0]}(111)=f_3+f_7,\ \  A^{[1]}(111)=A^{[0]}(011)+\omega^{2^2}A^{[0]}(111)=f_3-f_7
\end{eqnarray*}
  In pass 2,
   \begin{eqnarray*}
& &  A^{[2]}(000)=A^{[1]}(000)+\omega^{2\times(00)_2}A^{[1]}(010)=(f_0+f_4)+(f_2+f_6)\\
& &  A^{[2]}(001)=A^{[1]}(001)+\omega^{2\times(00)_2}A^{[1]}(011)=(f_1+f_5)+(f_3+f_7) \\
& &  A^{[2]}(010)=A^{[1]}(000)+\omega^{2\times(10)_2}A^{[1]}(010)=(f_0+f_4)-(f_2+f_6)\\
& &  A^{[2]}(011)=A^{[1]}(001)+\omega^{2\times(10)_2} A^{[1]}(011)=(f_1+f_5)-(f_3+f_7) \\
& &  A^{[2]}(100)=A^{[1]}(100)+\omega^{2\times(01)_2} A^{[1]}(110)=(f_0-f_4)+\omega^2(f_2-f_6)\\
& &  A^{[2]}(101)=A^{[1]}(101)+\omega^{2\times(01)_2}A^{[1]}(111)=(f_1-f_5)+\omega^2(f_3-f_7)\\
& &  A^{[2]}(110)=A^{[1]}(100)+\omega^{2\times(11)_2} A^{[1]}(110)=(f_0-f_4)-\omega^2(f_2-f_6) \\
& &  A^{[2]}(111)=A^{[1]}(101)+\omega^{2\times(11)_2}A^{[1]}(111)=(f_1-f_5)-\omega^2(f_3-f_7)
\end{eqnarray*}
 In pass 3,
 \begin{eqnarray*}
& &  A^{[3]}(000)=A^{[2]}(000)+\omega^{(000)_2}A^{[2]}(001)=[(f_0+f_4)+(f_2+f_6)]+[(f_1+f_5)+(f_3+f_7)]=f(1)\\
& &  A^{[3]}(001)=A^{[2]}(000)+\omega^{(100)_2}A^{[2]}(001)=[(f_0+f_4)+(f_2+f_6)]-[(f_1+f_5)+(f_3+f_7)]=f(\omega^4) \\
& &  A^{[3]}(010)=A^{[2]}(010)+\omega^{(010)_2}A^{[2]}(011)=[(f_0+f_4)-(f_2+f_6)]+\omega^2[(f_1+f_5)-(f_3+f_7)]=f(\omega^2)\\
& &  A^{[3]}(011)=A^{[2]}(010)+\omega^{(110)_2} A^{[2]}(011)=[(f_0+f_4)-(f_2+f_6)]-\omega^2[(f_1+f_5)-(f_3+f_7)]=f(\omega^6) \\
& &  A^{[3]}(100)=A^{[2]}(100)+\omega^{(001)_2} A^{[2]}(101)=[(f_0-f_4)+\omega^2(f_2-f_6)]+\omega[(f_1-f_5)+\omega^2(f_3-f_7)]=f(\omega)\\
& &  A^{[3]}(101)=A^{[2]}(100)+\omega^{(101)_2}A^{[2]}(101)=[(f_0-f_4)+\omega^2(f_2-f_6)]-\omega[(f_1-f_5)+\omega^2(f_3-f_7)]=f(\omega^5)\\
& &  A^{[3]}(110)=A^{[2]}(110)+\omega^{(011)_2} A^{[2]}(111)=[(f_0-f_4)-\omega^2(f_2-f_6)]+\omega^3[(f_1-f_5)-\omega^2(f_3-f_7)]=f(\omega^3)\\
& &  A^{[3]}(111)=A^{[2]}(110)+\omega^{(111)_2}A^{[2]}(111)=[(f_0-f_4)-\omega^2(f_2-f_6)]-\omega^3[(f_1-f_5)-\omega^2(f_3-f_7)]=f(\omega^7)
\end{eqnarray*} }
\emph{Drawback}: The description-2 almost veils the elegant idea behind FFT.

 \section{A heuristic description of FFT}

We now present an  explicit and heuristic description of FFT, which can be regarded as the combination of the above two descriptions.

The basic idea of the new description  is to recursively split polynomials and shrink coefficients. Concretely, given a polynomial $f(x)=f_0+f_1x+f_2x^2+\cdots+f_7x^7 \in \mathcal {C}[x] $ and a primitive 8th  root of unity $\omega=e^{2\pi i/8}$, we split it into two polynomials of  degree less than 4. Their coefficients can be individually obtained by shrinking the symmetric two coefficients of original polynomial.  See
the following working flow for details.

 \newpage

 \hspace*{-110mm}\begin{minipage}{\linewidth}
  \includegraphics[angle=90,height=21cm,width=25cm]{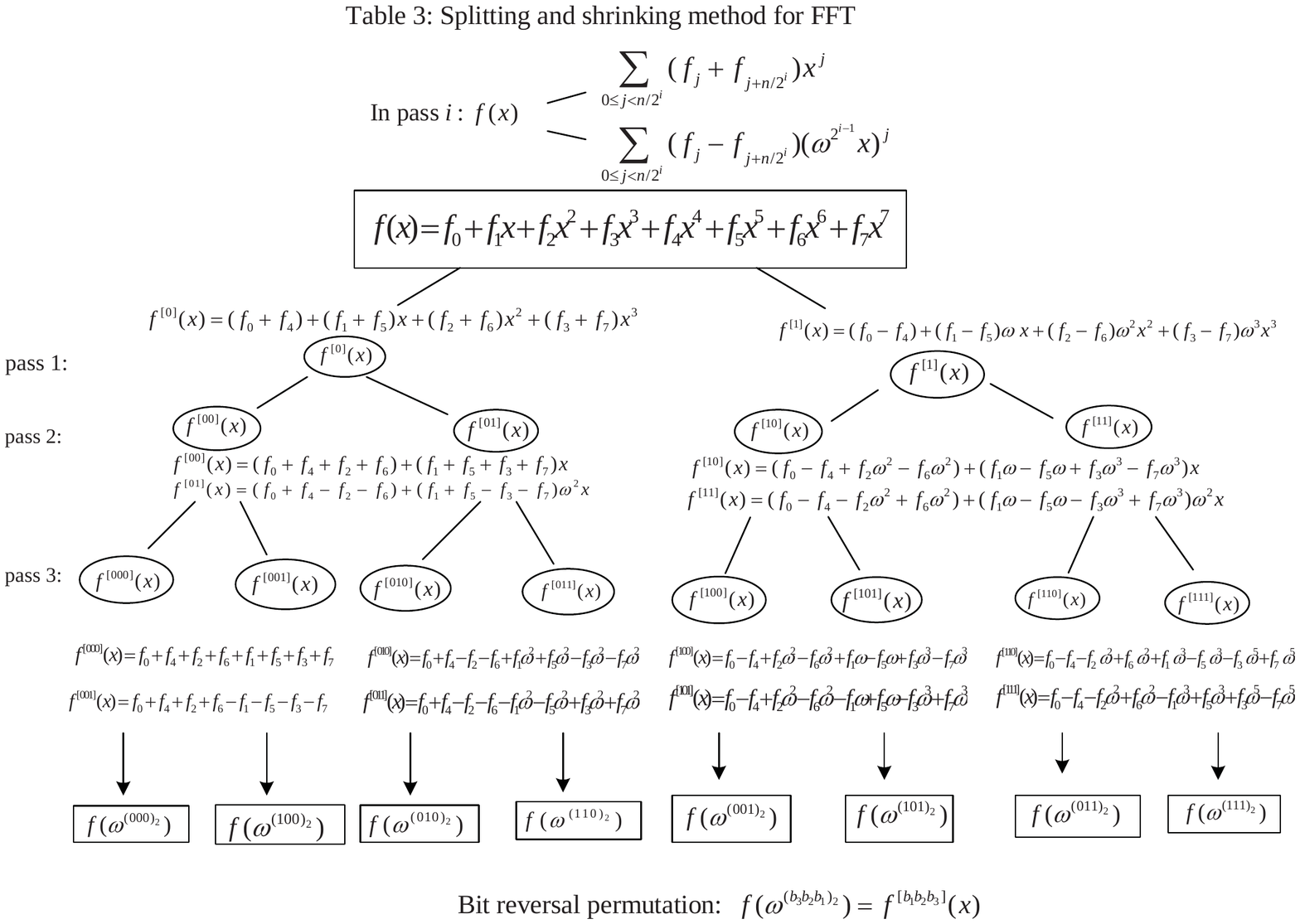}
 \end{minipage}

 \newpage

Refining the method, we have the following heuristic description for FFT:

  \begin{table}[!h]
\tabcolsep 0pt
\caption{Description-3}
\vspace*{-12pt}
\begin{center}
\def\temptablewidth{1\textwidth}
{\rule{\temptablewidth}{.7pt}}

\begin{tabular*}{\temptablewidth}{@{\extracolsep{\fill}}l}

Input: $n=2^k \in \mathcal {N}$ with $k\in \mathcal {N}, f(x)=\sum_{0\leq i<n}f_i\,x^i\in \mathcal {C}[x]$,
and \\  \hspace*{8mm} the powers $ \omega, \omega^2, \cdots, \omega^{n-1}$ of a primitive $n$th root of unity $\omega \in \mathcal {C}$.\\
Output:  $\mathrm{DFT}_\omega(f)=(f(1),f(\omega),\cdots,f(\omega^{n-1}))\in \mathcal {C}^n.$\\
\textbf{Pass 1.} Set $ f^{[b_1]}(x)= \left\{
\begin{array}{ll}
 \sum_{0\leq j<{n/2}}\left(f_j+f_{j+n/2}\right)x^j,  & b_1=0 \\
 \sum_{0\leq j<{n/2}}\left(f_j-f_{j+n/2}\right)(\omega\,x)^j, & b_1=1 \\
 \end{array}
 \right.
$ \\
\textbf{Pass 2.} For each $b_i\in \{0, 1\}, 1\leq i\leq 2$, compute\\
\hspace*{6mm} $ f^{[b_1b_2]}(x)= \left\{
\begin{array}{ll}
 \sum_{0\leq j<{n/2^2}}\left(f^{[b_1]}_j+f^{[b_1]}_{j+n/2^2}\right)x^j,  & b_2=0 \\
 \sum_{0\leq j<{n/2^2}}\left(f^{[b_1]}_j-f^{[b_1]}_{j+n/2^2}\right)(\omega^{2}\,x)^j, & b_2=1 \\
 \end{array}
 \right.
$\\
\hspace*{12mm} $\cdots$\\
\textbf{Pass k.} For each $b_i\in \{0, 1\}, 1\leq i\leq k$, compute\\
\hspace*{6mm} $ f^{[b_1b_2\cdots b_k]}(x)= \left\{
\begin{array}{ll}
 \sum_{0\leq j<n/2^{k}}\left(f^{[b_1b_2\cdots b_{k-1}]}_j+f^{[b_1b_2\cdots b_{k-1}]}_{j+n/2^{k}}\right)x^j,  & b_k=0 \\
  \sum_{0\leq j<n/2^{k}}\left(f^{[b_1b_2\cdots b_{k-1}]}_j-f^{[b_1b_2\cdots b_{k-1}]}_{j+n/2^{k}}\right)(\omega^{2^{k-1}}\,x)^j, & b_k=1 \\
 \end{array}
 \right. $ \\
\textbf{Reversal permutation.} For each $b_i\in \{0, 1\}, 1\leq i\leq k$,   set \\
\hspace*{6mm} $f(\omega^{(b_kb_{k-1}\cdots b_1)_2})= f^{[b_1\cdots b_{k-1}b_k]}(x)$ \\
\end{tabular*}
       {\rule{\temptablewidth}{.6pt}}
       \end{center}
       \end{table}

 It is easy to prove the correctness of above method  by induction and show that the algorithm requires $k\,n$ complex-number additions.
 Incidentally, we refer to Ref.[2] for an explicit pseudocode of FFT (ITERATIVE-FFT)

\section{Conclusion}

In this paper, we analyze two common descriptions of FFT and propose a heuristic presentation of complex-data FFT.
We think the new presentation is helpful for students and programmers to grasp the  method entirely and deeply.

 \emph{Acknowledgements}
 We thank the National Natural Science
Foundation of China (Project 60873227), and the Key Disciplines of
 Shanghai Municipality (S30104).

\end{document}